\documentclass[final,onefignum,onetabnum]{siamart171218}
\usepackage{ulem}

\usepackage{braket,amsfonts,amssymb,multirow}

\usepackage{array}

\usepackage[caption=false]{subfig}
\usepackage{pgfplots}

\newsiamthm{claim}{Claim}
\newsiamremark{remark}{Remark}
\newsiamremark{hypothesis}{Hypothesis}
\crefname{hypothesis}{Hypothesis}{Hypotheses}
\newsiamthm{ex}{Example}

\usepackage{algorithmic}

\usepackage{graphicx,epstopdf}

\Crefname{ALC@unique}{Line}{Lines}

\usepackage{amsopn}

\usepackage{xspace}
\usepackage{bold-extra}
\usepackage[most]{tcolorbox}

\colorlet{texcscolor}{blue!50!black}
\colorlet{texemcolor}{red!70!black}
\colorlet{texpreamble}{red!70!black}
\colorlet{codebackground}{black!25!white!25}

\lstdefinestyle{siamlatex}{%
  style=tcblatex,
  texcsstyle=*\color{texcscolor},
  texcsstyle=[2]\color{texemcolor},
  keywordstyle=[2]\color{texemcolor},
  moretexcs={cref,Cref,maketitle,mathcal,text,headers,email,url},
}

\tcbset{%
  colframe=black!75!white!75,
  coltitle=white,
  colback=codebackground, 
  colbacklower=white, 
  fonttitle=\bfseries,
  arc=0pt,outer arc=0pt,
  top=1pt,bottom=1pt,left=1mm,right=1mm,middle=1mm,boxsep=1mm,
  leftrule=0.3mm,rightrule=0.3mm,toprule=0.3mm,bottomrule=0.3mm,
  listing options={style=siamlatex}
}

\newtcblisting[use counter=example]{example}[2][]{%
  title={Example~\thetcbcounter: #2},#1}

\newtcbinputlisting[use counter=example]{\examplefile}[3][]{%
  title={Example~\thetcbcounter: #2},listing file={#3},#1}

\DeclareTotalTCBox{\code}{ v O{} }
{ 
  fontupper=\ttfamily\color{black},
  nobeforeafter,
  tcbox raise base,
  colback=codebackground,colframe=white,
  top=0pt,bottom=0pt,left=0mm,right=0mm,
  leftrule=0pt,rightrule=0pt,toprule=0mm,bottomrule=0mm,
  boxsep=0.5mm,
  #2}{#1}

\patchcmd\newpage{\vfil}{}{}{}
\begin{tcbverbatimwrite}{tmp_\jobname_header.tex}
\title{Topological properties of elastoplastic lattice spring models that determine terminal distributions of plastic deformations\thanks{Submitted to the editors DATE.
\funding{The authors were supported by NSF Grant CMMI-1916876.}}}

\author{Oleg Makarenkov\thanks{Department of Mathematical Sciences, University of Texas at Dallas, 75080 Richardson, USA (\email{makarenkov@utdallas.edu}, \email{josean.albelo-cortes@utdallas.edu}).}
\and Josean Albelo-Cortes\footnotemark[2]}

\end{tcbverbatimwrite}
\input{tmp_\jobname_header.tex}

\ifpdf
\hypersetup{ pdftitle={Guide to Using  SIAM'S \LaTeX\ Style} }
\fi

\def\qed{\hfill $\square$}

\begin{document}
\maketitle

\begin{tcbverbatimwrite}{tmp_\jobname_abstract.tex}
\begin{abstract}
A recent result by Gudoshnikov et al [SIAM J. Control Optim. 2022] ensures finite-time convergence of the stress-vector of an elastoplastic lattice spring model under assumption that the vector $g'(t)$ of the applied displacement controlled-loading lies strictly inside the normal cone to the associated polyhedral set (that depends on mechanical parameters of the springs). Determination of the terminal distribution of stresses has been thereby linked to a problem of spotting a face on the boundary of the polyhedral set where the normal cone contains vector $g'(t)$. In this paper the above-mentioned problem of spotting an eligible face is converted into a search for an eligible set of springs that have a certain topological property with respect to the entire graph (of springs). Specifically, we prove that eligible springs are those that keep non-zero lengths after their nodes are collapsed with the nodes of the displacement-controlled loading after a finite number of eligible displacements of the nodes of the graph. The proposed result allows to judge about possible distribution of plastic deformations in elastoplastic lattice spring models directly from the topology of the associated graph of springs. A benchmark example is provided. 
\end{abstract}


\end{tcbverbatimwrite}
\input{tmp_\jobname_abstract.tex}

\section{Introduction}

\begin{figure}[h]
\vskip-0.5cm
\centerline{\includegraphics[scale=0.55]{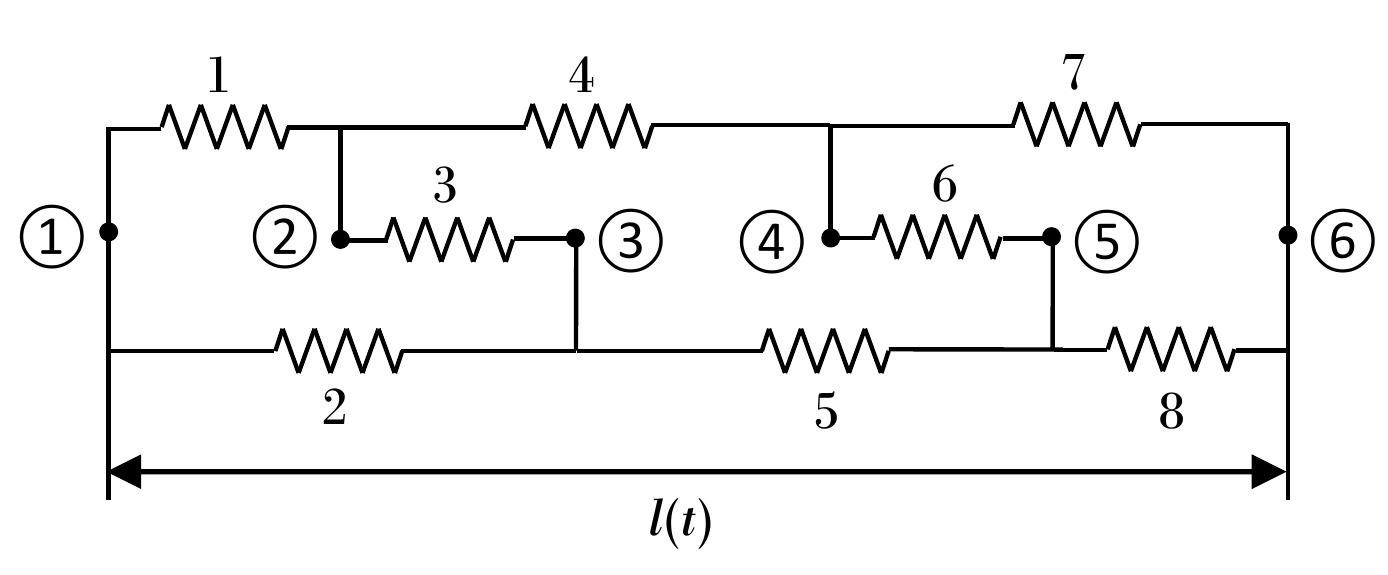}\qquad  \includegraphics[scale=0.55]{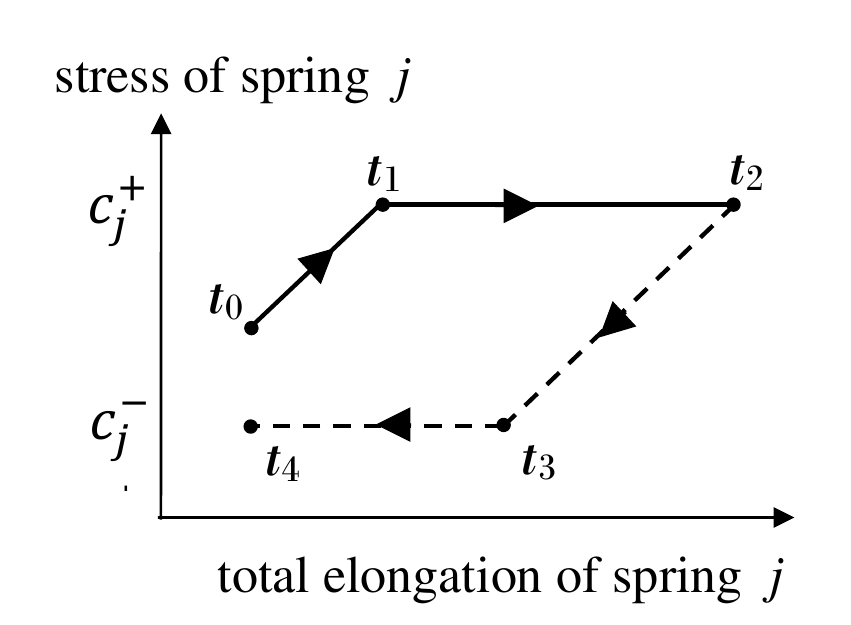}}
\caption{Left: A lattice spring model of 8 elastoplastic springs on 6 nodes with a unidirectional displacement-controlled loading $l(t)$. It is assumed that the $\xi$-axis is directed as in Fig.~\ref{figR}. Right: Illustration of elastoplastic properties of spring $j$ under tension (on interval $[t_0,t_2]$) and compression (on the interval $[t_2,t_4]$), and the role of the elastic limits $c_j^\pm.$ }\label{benchmark}
\end{figure}

\noindent A significant breakthrough in understanding of the response of elastoplastic lattice spring models \cite{JiaoPrior2,Jiao2,tissue,novel}  to displacement-controlled loading has recently been due to the development of the theory of sweeping processes \cite{GKMV,ESAIM,SICON,Moreau}, which links the asymptotic dynamics of the stress-vector $s(t)$ of the springs to an algebraic inclusion. 
Specifically, according to \cite{SICON}, 
when the displacement-controlled loading (see Fig.~\ref{benchmark}) is uni-directional, i.e.
$$
   l(t)=l_0+l_1t,
$$
the possible asymptotic time-evolutions of the stress-vector $s(t)$ are determined by solutions $
I_0\subset\{-1,1\}\times \overline{1,m}$ of 
\begin{equation}\label{blueinclusion}
\left(\begin{array}{c}
c_1\\
0_{m-n+1}\end{array}\right) \in {\rm cone}\left(
\left(\begin{array}{c}
  R^T \\
  (D^\perp)^T\end{array}\right) \left\{\alpha e_j:(\alpha,j)\in I_0\right\}\right)\quad{\rm with\ some}\ \ c_1>0,
\end{equation}
where $n\times m$-matrix $-D^T$ is the incidence matrix of the graph of springs and $R$ is the incidence vector of a path connecting the two endpoints of the displacement-controlled loading (see the beginning of the next Section). 
Indeed, if $I_0$ solves (\ref{blueinclusion}), then it is always possible to choose the elastic limits of springs $c_j^\alpha,$ $(\alpha,j)\in \{-1,1\}\times\overline{1,m}$ (see Fig.~\ref{benchmark}), in such a way that  $s(t)$ attains a constant value described by $I_0$, i.e. (\cite[Proposition 7.3]{SICON})
\begin{equation}\label{convergence}
 s(t)=s_*\ \ {\rm with}\ \  s_{*,j}=c_j^\alpha,\quad (\alpha,j)\in I_0,\quad \mbox {for all }t>0 \mbox{ sufficiently large,}
\end{equation}
where $s_{*,j}$ stays for the $j$-component of vector $s_*\in\mathbb{R}^m.$
In terms of mechanical properties of the lattice spring model, the equality $s_j(t)=c_j^\alpha$ means that spring $j$ reached its elastic limit and stays in plastic mode (capable to stretch plastically, if $\alpha=1$, and capable to compress plastically, if $\alpha=-1$), i.e. (\ref{blueinclusion}) is a key property in the determination of terminal distributions of plastic deformations.

\vskip0.2cm

\noindent The goal of this paper is to propose a purely topological rule that verifies validity of (\ref{blueinclusion}) for a given $I_0$ by looking directly at the graph of springs. We discover that if one collapses some of the nodes of the lattice spring model with one end of the displacement-control loading and collapses the remaining nodes with the other end of the displacement-controlled loading, then the indicies of springs that survive under such an operation can be taken as $I_0.$ For example, in the model of Fig.~\ref{benchmark}, collapsing the nodes as shown in Fig.~\ref{hhh}(a), springs 1, 3, 6, 8 keep non-zero length (survive) and springs 2, 4, 7 shrink to points (nullify), which we relate to the fact that a possible solution to (\ref{blueinclusion}) is $I_0=\{(\alpha_1,1),(\alpha_3,3),(\alpha_6,6),(\alpha_8,8)\}$ with suitable $\alpha_j\in\{-1,1\}$.  In other words, the paper proposes a rule to judge about possible distribution of plastic deformations in elastoplastic lattice spring models directly from the topology of the associated graph of springs.

\vskip0.2cm

\noindent The paper is organized as follows. In Section~\ref{genericity} we introduce some minimal background to explain (Theorem~\ref{maintheorem}) how (\ref{blueinclusion}) is related to plastic deformations in elastoplastic lattice spring models of the type of Fig.~\ref{benchmark}. 
The main result (Theorem~\ref{propmain}) about the determination of $I_0$ from the topology of the graph of springs is established in Section~\ref{mainresult}. An application of the main result to the lattice spring model of Fig.~\ref{benchmark} is discussed in Section~\ref{gameexample}. Conclusions section concludes the paper.

\section{Formulation of the finite-time stability theorem for elastoplastic lattice spring models}\label{genericity}

\noindent The goal of this section is to formulate the finite-time stability theorem of \cite{SICON} that concludes mechanical behavior (\ref{convergence}) from the algebraic inclusion (\ref{blueinclusion}). This will set a motivation behind (\ref{blueinclusion}) rigorously. 


\vskip0.2cm

\begin{figure}[h]\center
\vskip-0.0cm
\includegraphics[scale=0.47]{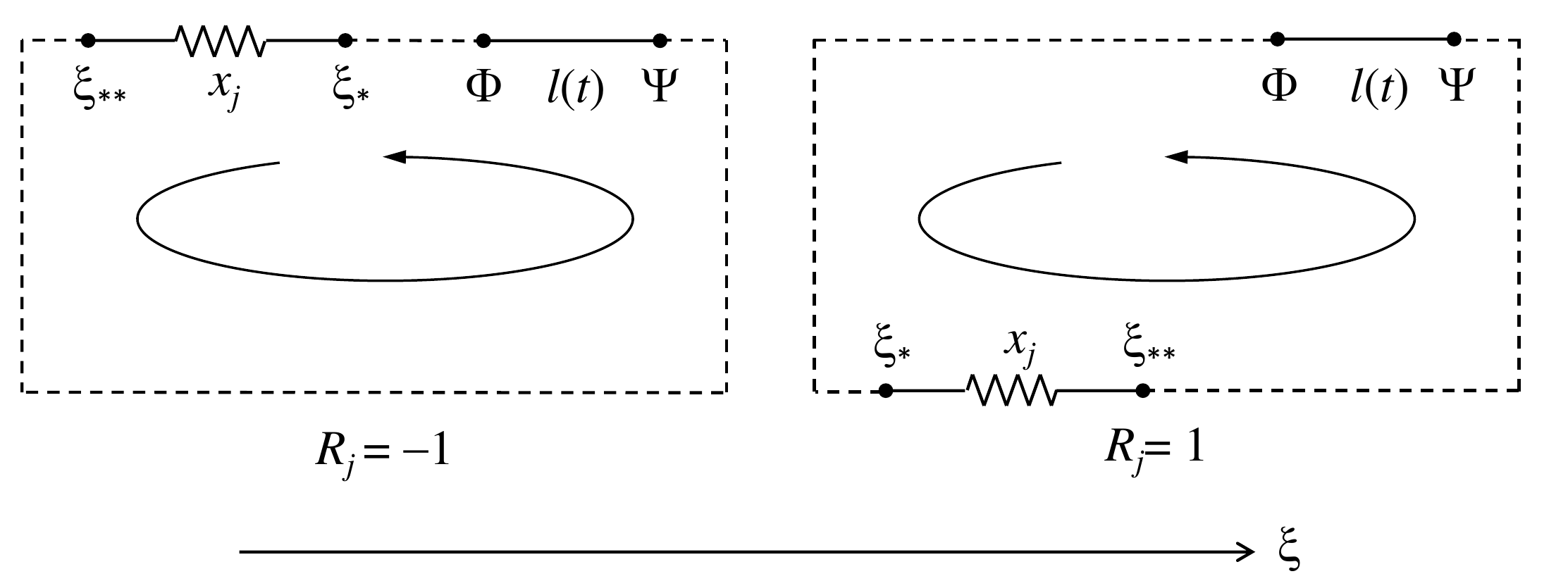}\vskip-0.0cm
\caption{ (follows \cite[Fig.~3]{PhysicaD}) Illustration of the signs of the components of the incidence vector $R\in\mathbb{R}^m.$ The dotted contour stays for the chain of the springs associated with the vector $R$, while each of the two figures discloses just one spring of the contour for illustration.
} \label{figR}
\end{figure}

\noindent The setup of this section follows  \cite{ESAIM,PhysicaD}. A lattice spring model of $m$ elastoplastic springs on $n$ nodes (whose coordinates are elements of $\mathbb{R}$) are connected according to a directed graph given by the $n\times m$ incidence matrix $-D^T$. Specifically, the $(i,j)$-element of matrix $D^T$ is 1 or $-1$ according to whether node $i$ is the terminus of spring $j$ or the source of spring $j$ \cite[Ch.~7]{Ivan22}. If none of these two cases takes place, then the $(i,j)$-element of matrix $D^T$ is 0. In other words, if $\xi\in\mathbb{R}^n$ is the vector of coordinates of the nodes of the springs, then $j$-th component of $D\xi$ is the length of spring $j$. This length is negative, when the coordinate of the terminus of spring $j$ is smaller than the coordinate of its source.
The Hooke's coefficients  $k_1, ..., k_m$ of the springs are arranged into a $k\times k$-matrix $K={\rm diag}\left\{k_1,...,k_m\right\}.$  
The displacement-controlled loading $l(t)$ is defined through a path of springs (with indicies, say, $j_1,...,j_s$) that connects one node of $l(t)$ (denoted by $\Phi$ in the Fig.~\ref{figR}) with its other node (denoted by $\Psi$). This path is described by a so-called {\it incidence vector} $R\in\mathbb{R}^m$ whose $j$-th component is 0 or $\pm1$ according to whether spring $j$ is a part of the path or not. If one follows the path beginning node $\Phi$ and heading towards node $\Psi$, and if $\xi_*$ and $\xi_{**}$ are two successive nodes on this way connected through spring $j$, then $j$-th component of the incident vector $R$ equals $-1$ or $1$ according to whether $\xi_*>\xi_{**}$ or $\xi_*<\xi_{**}$, see Fig.~\ref{figR}. Such an elastoplastic system will be referred to as $(D,K,C,R,l(t))$. A solution of lattice-spring model $(D,K,C,R,l(t))$ is an absolutely continuous function $\sigma(t)\in\mathbb{R}^m$ of the stresses of springs which verifies the following equations
\begin{align}
\text{Compatibility of springs' lengths:} & \  \ x\in D \mathbb{R}^n, \qquad 
 \label{D}\\
\text{Displacement-controlled loading:}  &\ \ R^Tx=l(t),\label{G}\\
\text{Additive plasticity of the lengths:} & \  \ x=K^{-1}\sigma+p,
 \label{addplast}\\
\text{Plastic deformation law:} 
&\ \ \dot p\in N_{[c_1^-,c_1^+]\times...\times[c_m^-,c_m^+]}(\sigma),\ \ {\rm a.\,e.\ on}\ [0,\infty), \label{formula2}\\
   \text{Static balance law:}&\ \  -D^T\sigma-D^TRr=0,\hskip1cm \label{eq5}
\end{align}
where $x(t)$, $p(t)$ and $r(t)$ are suitable absolutely continuous functions.

\vskip0.2cm 

\noindent Assuming that the graph of springs is connected, we have that (see Bapat \cite[Lemma 2.2]{Bapat})  $    {\rm rank}\hskip0.05cm D=n-1,$ i.e. matrix $D$ posses $n-1$ linearly independent columns. Therefore, there exists an $m\times(m-n+1)-$matrix $D^\perp$ that solves $  (D^\perp)^T D={\color{black}0_{(m-n+1)\times n}}$ and satisfies ${\rm rank}(D^\perp)=m-n+1. $

\begin{ex}  For the example of Fig.~\ref{benchmark}, the matrices $D$ and $D^\perp$ compute as 
$$
  D= \left(
\begin{array}{cccccc}
 -1 & 1 & 0 & 0 & 0 & 0 \\
 -1 & 0 & 1 & 0 & 0 & 0 \\
 0 & -1 & 1 & 0 & 0 & 0 \\
 0 & -1 & 0 & 1 & 0 & 0 \\
 0 & 0 & -1 & 0 & 1 & 0 \\
 0 & 0 & 0 & -1 & 1 & 0 \\
 0 & 0 & 0 & -1 & 0 & 1 \\
 0 & 0 & 0 & 0 & -1 & 1 \\
\end{array}
\right),\qquad D^\perp:={\rm Ker}\, D^T= \left(
\begin{array}{ccc}
 -1 & 1 & 1 \\
 1 & -1 & -1 \\
 0 & 0 & 1 \\
 -1 & 1 & 0 \\
 1 & -1 & 0 \\
 0 & 1 & 0 \\
 -1 & 0 & 0 \\
 1 & 0 & 0 \\
\end{array}
\right).
$$
Taking springs 1, 4, and 7 to create a path connecting the left and right ends of displacement-controlled loading $l(t)$, the incidence vector $R$ computes as
$$ R=\left(
\begin{array}{cccccccc}
 1 & 0 & 0 & 1 & 0 & 0 & 1 & 0 
\end{array}
\right)^T.
$$

\end{ex}

\vskip0.2cm


\begin{definition}\label{admissibilitydefinition} $I_0$ is admissible, if (\ref{blueinclusion}) holds.
\end{definition}

\begin{definition}\label{irreducibilitydefinition} An admissible $I_0$ is irreducible, if (\ref{blueinclusion}) fails for any $\tilde I_0\subset I_0$ such that $\tilde I_0\not=I_0$.
\end{definition}


\begin{theorem}\label{maintheorem} \mbox{\rm (\cite[Proposition~7.3]{SICON})} Assume that $I_0\subset \{-1,1\}\times\overline{1,m}$ is admissible and irreducible. Assume that 
\begin{equation}\label{assumptions}
  m>n,\qquad {\rm rank}(D^TR)=1.
\end{equation} 
Then, there exist $I_1,...,I_M\subset \{-1,1\}\times\overline{1,m}$, $M\in \mathbb{N}\cup\{0\},$ and $C_j^-<C_j^+$, such that if
$$
     c_j^-\le C_j^-,\quad C_j^+\le c_j^+,\quad  i\in\overline{1,M},\ (\alpha,j)\not\in I_0\cup I_1\cup ... \cup I_M,
$$
then 
every solution $\sigma(t)$ of the lattice spring model
$(D,K,C,R,l(t))$ satisfies (\ref{convergence}). 
\end{theorem}

\begin{remark} \label{remarkdim} According to \cite[Corollary~B.7]{SICON},
$
{\rm rank}\left(\begin{array}{c}
  R^T \\
  (D^\perp)^T\end{array}\right)
=\dim V$, which equals $m-n+2$. Therefore, the cardinality of an admissible and irreducible $I_0$ cannot exceed $m-n+2$. 
\end{remark}

\section{The main result}\label{mainresult}

\noindent Denote the indicies of $I_0$ as 
\begin{equation}\label{I0def}
  I_0=\{(\alpha_{j_1},j_1),...,(\alpha_{j_s},j_s)\}.
\end{equation}
Then inclusion (\ref{blueinclusion}) 
is equivalent to the existence of $\lambda_j\ge 0$ and $c_1>0$ such that the system of two relations holds:
\begin{eqnarray}
c_1 &=& R^T \lambda_{j_1}\alpha_{j_1}e_{j_1}+ \ldots+ R^T\lambda_{j_{s}}\alpha_{j_{s}}e_{j_{s}},\label{rel1}\\
0&=&(D^\perp)^T \left( \lambda_{j_1}\alpha_{j_1}e_{j_1}+ \ldots+ \lambda_{j_{s}}\alpha_{j_{s}}e_{j_{s}}\right).\label{rel2}
\end{eqnarray}

\begin{lemma}\label{lem0} If $I_0$ is admissible and irreducible, then $\lambda_{j_1},...,\lambda_{j_s}>0.$
\end{lemma}
{\bf Proof.} If $\lambda_{j_{s_*}}=0$, then (\ref{blueinclusion}) holds with $I_0$ replaced by $\tilde I_0=I_0\backslash\{(\alpha_{j_{s_*}},{j_{s_*}})\}.$\qed

\vskip0.2cm

\noindent 
By the definition of vector $R$, we have (see also Gudoshnikov et al \cite[formula (29)]{ESAIM})
\begin{equation}\label{PsiPhi}
  R^TD\xi=\xi_\Psi-\xi_\Phi,\quad \xi\in\mathbb{R}^n,
\end{equation}
where $\Phi,$ $\Psi$ are the indicies of the nodes of the displacement-controlled loading associated with $R$ (and chosen so that $\xi_\Phi<\xi_\Psi$), see the beginning of  Section~\ref{genericity}.

\begin{lemma}\label{lem1} (\ref{rel1})-(\ref{rel2}) holds if and only if there exist $\xi\in\mathbb{R}^n$ and $c_1>0$ 
such that
\begin{eqnarray}
   \xi_\Psi-\xi_\Phi&=&c_1, \label{relnewa} \\
    \lambda_{j_1}\alpha_{j_1}e_{j_1}+ \ldots+ \lambda_{j_{s}}\alpha_{j_{s}}e_{j_{s}} &=& D\xi. \label{relnewb} 
\end{eqnarray}
\end{lemma}

\noindent {\bf Proof.} Assume that (\ref{rel1})-(\ref{rel2}) hold. By (\ref{rel2}) and (see e.g. Friedberg et al. \cite[Exercise 17, p. 367]{Friedberg})
 \begin{equation}\label{weuse}
     {\rm Ker}\left((D^\perp)^T\right)=D\mathbb{R}^n,
 \end{equation}
 there exists $\xi\in\mathbb{R}^n$ such that (\ref{relnewb}) holds. By (\ref{rel1}) and (\ref{relnewb}),
\begin{equation}\label{romb}
   R^TD\xi=c_1,
\end{equation}
and (\ref{relnewa}) follows by combining (\ref{romb}) with (\ref{PsiPhi}).

\vskip0.2cm

\noindent Assume now that (\ref{relnewa})-(\ref{relnewb}) hold. Plugging (\ref{relnewa})-(\ref{relnewb})  to (\ref{romb}), one gets (\ref{rel1}) and (\ref{rel2}) follows from (\ref{weuse}).

\vskip0.2cm

\noindent The proof of the lemma is complete.\qed

\vskip0.3cm




\subsection{An informal formulation of the main result}
The main result of this paper is that (\ref{relnewa})-(\ref{relnewb}) holds, if each of the nodes of the lattice spring model $(D,K,C,R,l(t))$ can be collapsed with either node $\Phi$ or with node $\Psi$ in such a way that, after all nodes are collapsed with
nodes $\Phi$ or $\Psi$,
 springs $j_1,...,j_s$ keep non-zero lengths, and all the other springs shrink into points. One restriction that needs to be obeyed is that a node can be collapsed with $\Phi$ ($\Psi$) only if it is connected by a spring with a node whose  coordinate coincides with the coordinate of node $\Phi$ ($\Psi$) already.
For example, collapsing node 3 and then node 5 with node $1=\Phi$ and collapsing nodes 4 and then node 2 with node $6=\Psi$ in the lattice spring model of Fig.~\ref{benchmark},  
springs 1, 3, 6, 8 keep non-zero lengths and springs 2, 4, 7 shrink into points. Therefore, our result claims that, for the lattice spring model of Fig.~\ref{benchmark}, the property 
(\ref{relnewa})-(\ref{relnewb}) holds with 
$\{j_1,j_2,j_3,j_4\}=\{1,3,6,8\}.$ Furthermore, the corresponding values of $\alpha_{j_t}$ is $-1$ or $1$ according to whether the coordinates of the nodes of spring $j_t$ changed the order (on the $\xi$-axis) or not compared to the order before the nodes were collapsed.  As the diagram in Fig.~\ref{hhh}(a) illustrates, for the example under consideration, the order of nodes is preserved for springs 1 and 8, and the order is reversed for springs 3 and 6, that is why we conclude $\{\alpha_{j_1},\alpha_{j_2},\alpha_{j_3},\alpha_{j_4}\}=\{1,-1,-1,1\}$, i.e. (\ref{relnewa})-(\ref{relnewb}) holds with $I_0=\{(1,1),(-1,3),(-1,6),(1,8)\}.$  

\vskip0.2cm

\noindent In other words, taking into account Theorem~\ref{maintheorem}, the main result of the paper is that we can provide a candidate terminal distribution of plastic deformations in lattice spring model $(D,K,C,R,l(t))$ by collapsing the nodes of our choice with node $\Phi$ and by collapsing the remaining nodes with node $\Psi,$ where $\Phi$ and $\Psi$ are the indices of the two nodes of the displacement-controlled loading $l(t)$. 

\subsection{A formal formulation of the main result and a proof}
As announced in the previous section, we are going to show that 
in order for $\xi\in\mathbb{R}^n$ to be a solution of  
(\ref{relnewa})-(\ref{relnewb}),  $\xi$ must nullify the springs with the indices $\overline{1,m}\backslash\{j_1,...,j_s\}$ when the coordinates of the nodes of the lattice spring model $(D,K,C,R,l(t))$ are set to $\xi.$ 

\vskip0.2cm

\noindent 
Recall that when the coordinates of the nodes of the 
the lattice spring model $(D,K,C,R,l(t))$ are set to $\xi,$ 
the {\it length of spring $j$} 
is 
the $j$th component of the vector $D\xi.$
 Rephrasing the corresponding definition from graph theory, an $(i_1i_q)-walk$ is a sequence of nodes $i_1,...,i_q$ such that nodes $i_s$ and $i_{s+1}$ are connected by a spring ({\it adjacent}) for every $s\in\overline{1,q-1}$.

\begin{lemma}\label{necessarylemma} If $I_0$ is admissible and irreducible, then 
\begin{equation}\label{twostates}
   \xi_i\in\{\xi_\Phi,\xi_\Psi\}, \quad\mbox{for any}\ i\in\overline{1,n},
\end{equation}
and 
\begin{equation}\label{twostates2}
  \begin{array}{l}
\mbox{for any}\ i\in\overline{1,n},\\
      \mbox{if }\xi_i=\xi_\Phi\ { and}\ i\not=\Phi,\mbox{ then there exists an }(i\Phi)-walk 
\\
whose \ all \ nodes\ j\ satisfy \ \xi_j=\xi_\Phi,\\
      and\ same\ property\ holds\ by \ replacing\ \Phi\ by\ \Psi,
  \end{array}
\end{equation}
for any 
$\xi\in\mathbb{R}^n$ and any $\lambda_{j_1},...,\lambda_{j_{s}}\ge 0$ that satisfy (\ref{relnewa})-(\ref{relnewb}). 
\end{lemma}
  
\noindent {\bf Proof.} {\it Proof of (\ref{twostates}).} Let $I_0$ be irreducible and assume, by contradiction, that   
there exist $\xi\in\mathbb{R}^n$, $\lambda_{j_1},...,\lambda_{j_{s}}\ge 0$, and $i_*\in\overline{1,n},$ such that 
\begin{equation}\label{twostates-}
   \xi_{i_*}\not\in\{\xi_\Phi,\xi_\Psi\},
\end{equation} 
Observe that by connectivity of the graph at least one of the following two options  always holds:
\begin{itemize}
\item[A:] among all $i\in\overline{1,n}$ such that $\xi_i\not\in \{\xi_\Phi,\xi_\Psi\}$ there exists such $i\in\overline{1,n}$ for which $i$ and $\Phi$ are adjacent;
\item[B:] among all $i\in\overline{1,n}$ such that $\xi_i\not\in \{\xi_\Phi,\xi_\Psi\}$ there exists such $i\in\overline{1,n}$ for which $i$ and $\Psi$ are adjacent.
\end{itemize}

\noindent We will conduct the proof assuming that option A takes place. The case of option~B can be addressed by analogy.

\vskip0.2cm

\noindent We will now describe an iteration that transforms $\xi$ to some $\tilde \xi$. We then rename $\tilde \xi$ to $\xi$ and run the iteration again. This process will be repeated (a finite number of times) until the iteration no longer changes $\xi.$

\vskip0.2cm

\noindent {\bf Iteration of the process.}   Put
$$
   i_{**}=\underset{i\in\overline{1,n},\ i\not=\Phi}{\rm argmin}\{\xi_i\}.
$$
We have $\xi_{i_{**}}>\xi_\Phi$ by (\ref{twostates-}). 
\begin{equation}\label{action}
\begin{array}{l}
\mbox{Define }\tilde \xi\mbox{ by setting all components of }\xi\mbox{ which equal }\xi_{i_{**}}\\ \mbox{(including the component }\xi_{i_{**}}\mbox{ itself) to }\xi_\Phi. 
\end{array}
\end{equation}
Let $j_1,...,j_s$ be as given by (\ref{I0def}). We claim that $\tilde\xi$ satisfies (\ref{relnewb}) with some $\tilde\lambda_{j_1},...,\tilde\lambda_{j_{s}}\ge 0.$ Indeed, let $k\in\overline{1,m}$ be an arbitrary line of equality (\ref{relnewb}). We have 3 options to investigate how the length of spring $k$ changes 
under the replacement of $\xi$ by $\tilde\xi$:
\begin{itemize}
\item[1)] if exactly one of the nodes of spring $k$ is at position $\xi_{i_{**}}$, then the length of spring $k$ either resets to zero (if another node of spring $k$ is $\xi_\Phi$) or increases its nonzero absolute value (if another node of spring is not $\xi_\Phi$);
\item[2)] if both nodes of spring $k$ are at position $\xi_{i_{**}}$, then the length of spring $k$ doesn't change;
\item[3)] if none of the nodes of spring $k$ are at position $\xi_{i_{**}}$, then the length of spring $k$ doesn't change.
\end{itemize}
The conclusion of options 1)-3) is that only springs of nonzero length are capable to change length under the action of (\ref{action}) and if the length changes, the sign of the length doesn't change.  Therefore,
\begin{equation}\label{tildexi}
    \tilde\lambda_{j_1}\alpha_{j_1}e_{j_1}+ \ldots+ \tilde\lambda_{j_{s}}\alpha_{j_{s}}e_{j_{s}} = D\tilde\xi,
\end{equation}
where $\tilde\lambda_{j_1},...,\tilde\lambda_{j_s}\ge 0$.

\vskip0.2cm

\noindent 
Exit the iteration process if 
\begin{equation}\label{exit}
\tilde\lambda_{j_t}= 0, \mbox{ for some }t\in\overline{1,s}.
\end{equation} 
This event occurs in at most $n-2$ iterations by assumption of Option~A. Let $\xi:=\tilde \xi$ and iterate the process one more time, if it doesn't exit.

\vskip0.2cm 
\noindent According to Lemma~\ref{lem1} property (\ref{exit}) ensures that $\tilde I_0=I_0\backslash\{\lambda_{j_s},j_s\}$ is admissible. This means that we arrived to a contradiction with irreducibility of $I_0$ and proved (\ref{twostates}) for all $\xi\in\mathbb{R}^n$ and any $\lambda_{j_1},...,\lambda_{j_{s}}\ge 0$ that satisfy (\ref{relnewa})-(\ref{relnewb}).

\vskip0.2cm

\noindent Proof of (\ref{twostates2}).        Let $I_0$ be irreducible and assume, by contradiction, that   
there exist $\xi\in\mathbb{R}^n$, $\lambda_{j_1},...,\lambda_{j_s}\ge 0$, and $i_*\in\overline{1,n},$ such that 
$
\xi_{i_*}=\xi_\Phi$  and $i_*\not=\Phi$ while nodes $i_*$ and $\Phi$ are not connected      by  a  walk  whose  any  node $j$ satisfies $\xi_j=\xi_\Phi$. Let $J$ be the set of all nodes that are connected to node $i_*$ through walks with nodes at position $\xi_\Phi$ only. Since the graph of springs is connected, then there exists a spring (say, spring $k$) that connects $J$ to a node with coordinate $\xi_\Psi.$ Therefore, spring $k$ has nonzero length. Therefore, $k\in\{j_1,...,j_s\}.$ Create $\tilde \xi$ by moving nodes $J$ of $\xi$ from position $\xi_\Phi$ to position $\xi_\Psi.$ None of springs increase the length under this move because, by the definition of $J$, there are no springs that connect $J$ with $\xi_\Phi$, while springs that connect $J$ with $\xi_\Psi$ only decrease to $0$. In particular, spring $k$ gets zero length, which means that $\tilde\xi$ satisfies (\ref{tildexi}) with $\tilde \lambda_k=0$, contradicting (consult Lemma~\ref{lem1} again) irreducibility of $I_0$. 
 
 \vskip0.2cm
 
\noindent The proof of the lemma is complete.\qed

\begin{lemma}\label{second} Let $\xi\in\mathbb{R}^n$ satisfy (\ref{relnewa}) and (\ref{twostates}). Let $\{j_1,...,j_s\}\subset\overline{1,m}$ be all indices $j\in\overline{1,m}$ such that 
$$
  [D\xi]_{j}\not=0,
$$
so that (\ref{relnewb}) holds with 
$$
   \lambda_{j_k}=|[D\xi]_{j_k}|,\ \alpha_{j_k}={\rm sign}([D\xi]_{j_k}), \quad k\in\overline{1,s},
$$
i.e. $I_0$ is admissible.
If $\xi$ satisfies (\ref{twostates2}), then $I_0=\{(\alpha_{j_1},j_1),...,(\alpha_{j_s},j_s)\}$ is irreducible.
\end{lemma}

\noindent {\bf Proof.} Assume that $I_0$ is reducible and let $\tilde I_0\subset I_0$ be irreducible. Without loss of generality we can take
$$
    \tilde I_0=\{(\alpha_1,j_1),...,(\alpha_{j_t},j_t)\},
$$
where $t<s.$
By Lemma~\ref{lem1} there exists $\tilde\xi$ such that 
\begin{equation}\label{st}
   \tilde \lambda_{j_1}\alpha_{j_1}e_{j_1}+\ldots+\tilde \lambda_{j_t}\alpha_{j_t}e_{j_t}=D\tilde\xi\quad {\rm with}\ \tilde\lambda_{j_1},...,\tilde\lambda_{j_t}\ge 0.
\end{equation}
Since $\tilde I_0$ is irreducible, Lemma~\ref{necessarylemma} implies that $\tilde\xi_i\in\{\tilde\xi_\Phi,\tilde\xi_\Psi\},$ and by Lemma~\ref{lem1}, $\tilde \xi_\Psi-\tilde\xi_\Phi=c_1.$ We have
$$
  \tilde\xi_i+\xi_\Phi-\tilde\xi_\Phi\in\{\xi_\Phi,\tilde\xi_\Psi-\tilde\xi_\Phi+\xi_\Phi\}=\{\xi_\Phi,\xi_\Psi\},\quad i\in\overline{1,n}.
$$
Therefore, without loss of generality, we can assume that 
\begin{equation}\label{stasta}
   \tilde \xi_i\in\{\xi_\Phi,\xi_\Psi\},\quad i\in\overline{1,n}
\end{equation}
(i.e. we can subtract $\xi_\Phi-\tilde\xi_\Phi$ from all nodes of $\tilde\xi$ and denote by $\tilde\xi$ the vector obtained).

\vskip0.2cm

\noindent Equality (\ref{st}) and Lemma~\ref{lem0} imply that
$$
  [D\tilde\xi]_k=0\quad{\rm and}\quad [D\xi]_k\not=0,\quad k\in\overline{j_{t+1},...,j_s}.
$$ 
Therefore, moving the nodes from position $\xi$ to position $\tilde\xi$ changes the lengths of springs $j_{t+1},...,j_s$ from nonzero to zero values. Therefore, at least one node (say, node $i$) moves from position $\xi_\Phi$ to position $\xi_\Psi$ or from position $\xi_\Psi$ to position $\xi_\Phi.$ Assume that the former case takes place (the latter case can be considered by analogy). By condition (\ref{twostates2}) there is a path between node $i$ and node $\Phi$ of springs of zero lengths before the move. Therefore, one of the springs of this path must take nonzero value after the move. As a consequence, the left-hand-side of (\ref{st}) will get a spring of nonzero length that is not present in the left-hand-side of (\ref{relnewb}) (note that by the definition of $\{j_1,...,j_s\}$ in the formulation of the lemma, each $j_k$ corresponds to a spring on nonzero length in the left-hand-side of (\ref{relnewb})). Contradiction with the fact that $\tilde I_0\subset I_0.$ The proof of the lemma is complete.\qed

\begin{theorem} \label{propmain} Let $\xi\in\mathbb{R}^n$ be an arbitrary chosen vector of the coordinates of the nodes of the lattice spring model $(D,K,C,R,l(t))$ such that $\xi_i\not=\xi_\Phi,$ $i\not=\Phi$, and $\xi_i\not=\xi_\Psi,$ $i\not=\Psi$. Assume that the direction of the $\xi$-axis is chosen so that $\xi_\Phi<\xi_\Psi$. If a finite sequence of the following transformations 
\begin{equation}\label{trans}
\begin{array}{l}\mbox{pick any two nodes }i\mbox{ and }j\mbox{ adjoint one another and such that }\\
\xi_i\not\in\{\xi_\Phi,\xi_\Psi\}\mbox{ and }\xi_j\in\{\xi_\Phi,\xi_\Psi\};  
\mbox{ set }\xi_i\mbox{ to }\xi_j,\mbox{ i.e. put }\xi_i:=\xi_j,
\end{array}
\end{equation}
brings the original $\xi$ to a new $\xi$ that satisfies condition (\ref{twostates}), then the new $\xi$ satisfies condition (\ref{twostates2}). Furthermore, if $I_0=\{
(\alpha_{j_1},j_1),...,(\alpha_{j_s},j_s)\}$ is constructed for the new $\xi$ according to Lemma~\ref{second}, then $I_0$ is admissible and irreducible. 
\end{theorem}



\noindent {\bf Proof.} Since by the assumption of the theorem $\xi_i\not=\xi_\Phi$, $i\not=\Phi$, and $\xi_i\not=\xi_\Phi$, $i\not=\Phi$, the  given $\xi$ satisfies (\ref{twostates2}). We will proceed by induction.

\vskip0.2cm

\noindent Assume that $\xi$ satisfies (\ref{twostates2}) and let $\tilde\xi$ be obtained from $\xi$ by applying (\ref{trans}) with some $i_*,j_*\in\overline{1,n}.$ Since $\xi$ satisfies (\ref{twostates2}), then there exists an $(j_*\Phi)$-walk $j_*,i_1,...,i_q,\Phi$ such that
$$
  \xi_{j_*}=\xi_{i_1}=...=\xi_{i_q}=\xi_\Phi.
$$
Since (\ref{trans}) changes the value of one element $\xi_{i_*}\not\in\{\xi_\Phi,\xi_\Psi\},$ we have
$$
   \tilde \xi_{i_{j}}=\xi_{i_{j}},\ \ \tilde\xi_{i_1}=\xi_{i_1},\ \ ...\ \ \tilde\xi_{i_q}=\xi_{i_q},\ \ \tilde\xi_\Phi=\xi_\Phi.
$$
Therefore, $i_*,j_*,i_1,...,i_q,\Phi$ is a $(i_*\Phi)$-walk satisfying 
$$
  \tilde\xi_{i_*}=\tilde\xi_{j_*}=\tilde\xi_{i_1}=...=\tilde\xi_{i_q}=\tilde\xi_\Phi,
$$
i.e. $\tilde \xi$ satisfies (\ref{twostates2}) for $i=i_*$. 

\vskip0.2cm

\noindent Let now $i_{**}$ be any node of $\overline{1,n}$ other than $i_*$ and such that either (a) $\xi_{i_{**}}=\xi_\Phi$ and $i_{**}\not=\Phi$ or  (b) $\xi_{i_{**}}=\xi_\Psi$ and $i_{**}\not=\Psi$. Assume that case (b) takes place. Case (a) can be dealt with be analogy. By (\ref{twostates2}), there exists an $(i_{**}\Psi)$-walk $i_{**},i_1,...,i_q,\Psi$ such that
$$
  \xi_{i_{**}}=\xi_{i_1}=...=\xi_{i_q}=\xi_\Psi.
$$
As above, since (\ref{trans}) changes the value of one element $\xi_{i_*}\not\in\{\xi_\Phi,\xi_\Psi\},$ we conclude that 
$$
  \tilde\xi_{i_{**}}=\tilde\xi_{i_1}=...=\tilde\xi_{i_q}=\tilde\xi_\Psi.
$$
The proof of the fact that $\tilde \xi$ satisfies (\ref{twostates2}) is now complete.

\vskip0.2cm

\noindent  Assuming additionally, that repeated applications of transformation (\ref{trans}) brings the initial $\xi$ to a $\tilde \xi$ satisfying (\ref{twostates}), we have that $\tilde \xi$ satisfies both (\ref{twostates}) and (\ref{twostates2})

\vskip0.2cm

\noindent The requirement (\ref{relnewa}) is satisfied by setting 
$$
   c_1=\tilde\xi_\Psi-\tilde\xi_\Phi,
$$
so that the conclusion follows by applying Lemma~\ref{second}.

\vskip0.2cm

\noindent The proof of the theorem is complete.\qed

\section{Application to the benchmark example}\label{gameexample}

\begin{figure}
\centerline{\includegraphics[scale=0.55]{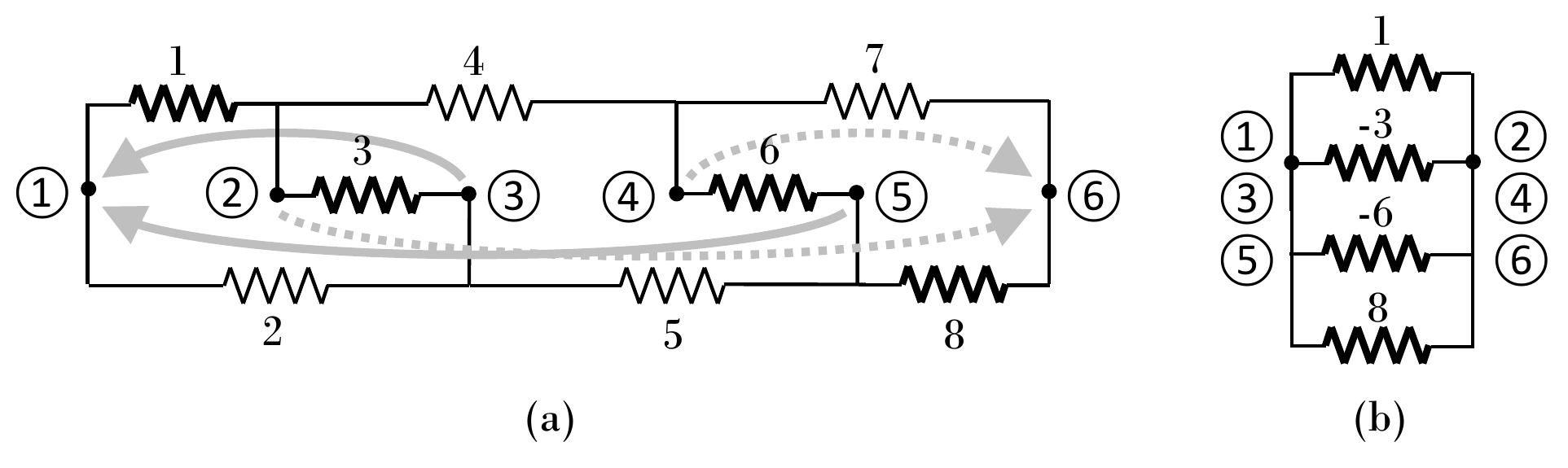}}
\centerline{\includegraphics[scale=0.55]{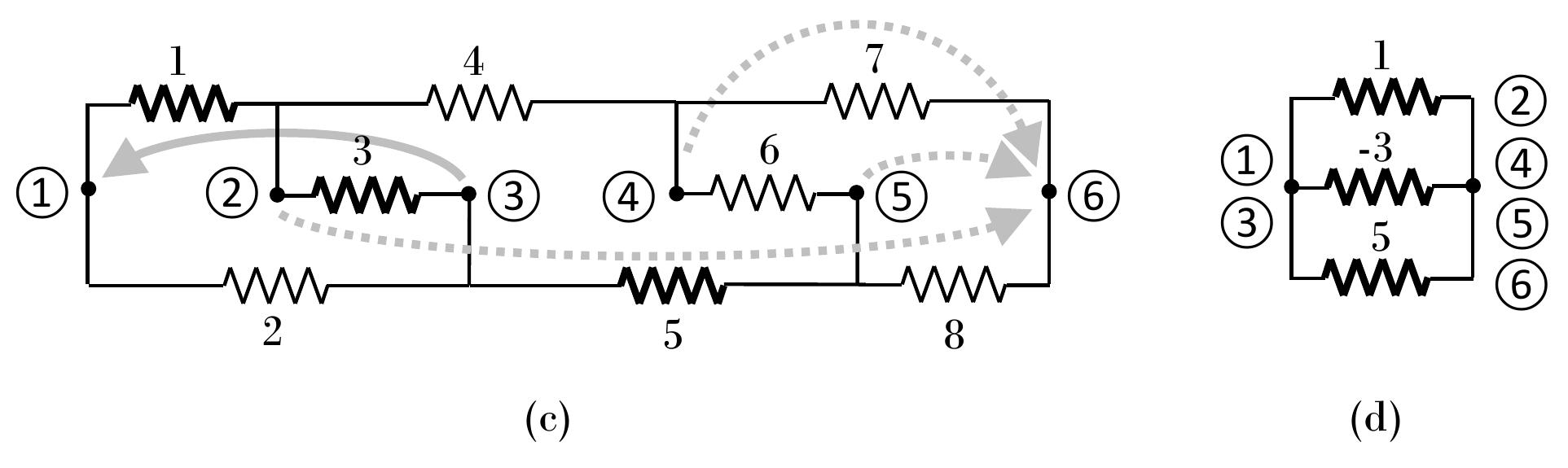}}
\centerline{\includegraphics[scale=0.55]{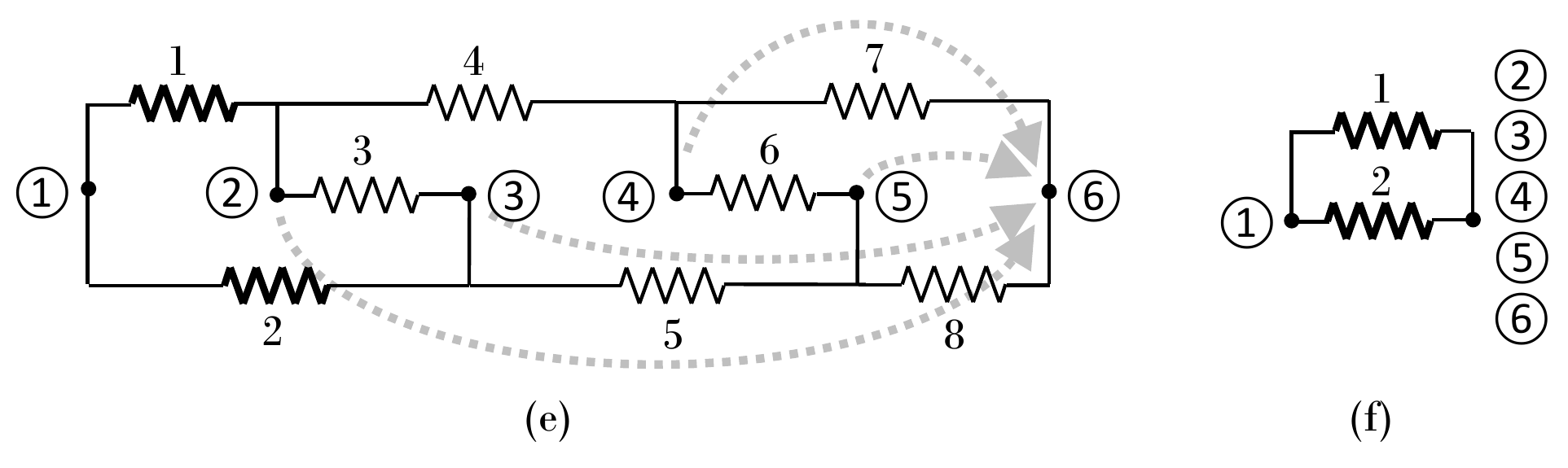}}
\caption{Left-hand figures: Subsequent displacements of the nodes of the lattice spring model of Fig.~\ref{benchmark} 
that satisfy the requirement (\ref{trans}); Right-hand figures: The configuration of the lattice spring model after all nodes are collapsed according to the rules from the corresponding left-hand figure. It is assumed that the $\xi$-axis is directed as in Fig.~\ref{figR}.}
\label{hhh}
\end{figure}

\noindent According to Remark~\ref{remarkdim}, the maximal cardinality of an admissible and irreducible $I_0$ 
for the example of Fig.~\ref{benchmark} is $m-n+2=4$. Therefore, we will now use Theorem~\ref{propmain} in order to spot 4 springs $j_1,j_2,j_3,j_4$ that do not shrink to a point under 
transformations (\ref{trans}) but get same coordinates of the associated nodes in the sense of the property (\ref{twostates}). Indeed, let us consider the following sequence of transformations (\ref{trans}) 
\begin{equation}\label{transexample}
  \xi_3:=\xi_1,\ \ \xi_5:= \xi_3,\ \ \xi_4:= \xi_6,\ \ \xi_2:=\xi_4,
\end{equation}
where the shortcut $\xi_i:=\xi_j$ stays for "set $\xi_i$ to $\xi_j$". One can check that $\xi_3:=\xi_1$ satisfies (\ref{transexample}) because nodes $1$ and $3$ are adjacent. While $\xi_5:=\xi_3$ wouldn't satisfy (\ref{transexample}) initially (i.e. before $\xi_3:=\xi_1$ is executed) since $\xi_3\not\in\{\xi_1,\xi_6\}$ initially, the property $\xi_3\in\{\xi_1,\xi_6\}$ holds after $\xi_3:=\xi_1$ is executed, so that $\xi_5:= \xi_3$ is eligible as being applied after $\xi_3:=\xi_1$. Analogously $\xi_2:=\xi_4$ is eligible because it is applied after $\xi_4:= \xi_6.$
The 4 transformations of (\ref{transexample}) are illustrated at Fig.~\ref{hhh}(a). The nodes for the final vector $\xi$ are drawn at Fig.~\ref{hhh}(b). In Fig.~\ref{hhh}(b) we added a negative sign to the index $j$ of a spring when transformation (\ref{transexample}) changed $[D\xi]_j$ from positive to negative (note, $[D\xi]_j>0$, $j\in\overline{1,8}$, for the original $\xi$ that corresponds to the diagram of Fig.~\ref{benchmark}). These springs are spring~3 and spring~6. Therefore, applying Theorem~\ref{propmain} with transformations (\ref{transexample}), we conclude that an admissible and irreducible $I_0$ is given by
\begin{equation}\label{I0-1}
  I_0=\{(1,1),(-1,3),(-1,6),(1,8)\}.
\end{equation}

\vskip0.2cm

\noindent Using Theorem~\ref{propmain} to discover $I_0$ with $|I_0|<4$ is a simpler task. For example, applying Theorem~\ref{propmain} with the sequence of transformations 
$$
   \xi_3:=\xi_1,\ \ \xi_5:=\xi_6;\ \ \xi_4:=\xi_5,\ \ \xi_2:=\xi_4,
$$
leads to the following admissible and irreducible $I_0$:
\begin{equation}\label{I0-2}
  I_0:=\{(1,1),(-1,3),(1,5)\},
\end{equation}
as illustrated in Fig.~\ref{hhh}(c)-(d). An $I_0$ of cardinality 2 
\begin{equation}\label{I0-3}
  I_0:=\{(1,1),(1,2)\},
\end{equation}
can be obtained over the following sequence of transformations
$$
   \xi_5:=\xi_6,\ \ \xi_4:=\xi_6;\ \ \xi_3:=\xi_5,\ \ \xi_2:=\xi_4,
$$
Fig.~\ref{hhh}(e)-(f).

\vskip0.2cm

\noindent To summarize, by using just topological observations of Fig.~\ref{hhh}, we were able to conclude that each of the sets (\ref{I0-1})-(\ref{I0-3}) defines a realizable distribution of plastic deformations in the lattice spring model of Fig.~\ref{benchmark} in the sense of Theorem~\ref{maintheorem}.

\section{Conclusions}\label{concl}

\noindent We proposed a topological rule to predict possible distributions of plastic deformations in elastoplastic lattice spring models. The approach doesn't require any algebraic computations and is based on a sequence of displacements of the nodes according to rule (\ref{trans}). Different eligible sequences of displacements of nodes correspond to different possible distributions of plastic deformations. 

\vskip0.2cm

\noindent The proposed topological rule relates the verification of algebraic inclusion (\ref{blueinclusion}) to the problem of bi-partitioning of connected graphs (our work in progress \cite{Galveston}). Though we know that the maximal cardinality of $I_0$ is $\dim V=m-n+2$, we don't know whether a solution $I_0$ of (\ref{blueinclusion}) of cardinality $m-n+2$ always exists. And if $I_0$ of cardinality $m-n+2$ exists, how many different $I_0$ with $|I_0|=m-n+2$ does a given lattice spring model $(D,K,C,R,l(t))$ allow? In the example of Fig.~\ref{benchmark}, one can easily find $|I_0|=4,$ $|I_0|=3$, and $|I_0|=2$ other than those considered in Fig.~\ref{hhh}. If there is a lattice spring model for which $|I_0|=m-n+2$ satisfying (\ref{blueinclusion}) doesn't exist, what is the maximal cardinality of $I_0$ that is capable to solve (\ref{blueinclusion})? We expect that this paper will attract interest of experts in graph theory.

\vskip0.2cm

\noindent On the other hand, answers to the above questions will help engineers to 
design networks of elastoplastic springs that allow most uniform distribution of plastic deformations across the material which minimizes the risk of crack initialization \cite{Blechman}.

\vskip0.2cm 

\noindent Inclusion (\ref{blueinclusion}) has been obtained in \cite{SICON} as a condition for a certain differential inclusion with moving constraint (sweeping process) to converge to a terminal solution in finite time. Specifically, the work \cite{SICON} established an abstract condition for finite-time stability of a sweeping process with a polyhedral moving constraint. We were able to reduce this abstract condition to the form (\ref{blueinclusion}) in \cite{SICON} (that we are now able to verify over the topological rule of Theorem~\ref{propmain}) only in the case where the sweeping process comes from a lattice spring model (with one-dimensional nodes). It is unknown whether finite-time stability of other types of sweeping processes with polyhedral moving constraint (e.g. those modeling dry friction systems 
\cite{Adly,Gidoni}) can be verified over some kind of adaptation of Theorem~\ref{propmain}. We hope that this line of questions will attract the interest of experts in stability of sweeping processes.




%


\section*{Acknowledgments} We thank Ivan Gudoshnikov (Czech Academy of Sciences) who read the manuscript and provided several important remarks. 






\end{document}